\documentstyle[amssymb,amsmath]{article}

\def\OE{\O}

\begin{document}

\renewcommand{\le}{\leqslant}
\renewcommand{\ge}{\geqslant}
\renewcommand{\epsilon}{\varepsilon}

\def\OE{\O}

\def\R {{\Bbb R }}
\def\C {{\Bbb C }}
\def\M{{\rm M}}
\def\OM{\overline \M}
\def\NN{{\cal N}}
\def\SS{{\sc S}}
\def\F{{\cal F}}
\def\H{{\cal H}}

\def\GL{{\rm GL}}
\def\U{{\rm U}}
\def\O{{\rm O}}
\def\Sp{{\rm Sp}}
\def\const{{\rm const}}
\def\B{{\rm B}}
\def\SO{{\rm SO}}
\def\OB{\overline {\B}}
\def\L{L^1_{\rm loc}}
\def\th{{\rm th\,\,}}

\newcounter{sec}
 \renewcommand{\theequation}{\arabic{sec}.\arabic{equation}}
\newcounter{fact}
\newcommand{\fact}[1]{\addtocounter{fact}{1} {\sc #1 \arabic{sec}.\arabic{fact}.}}

%{\sc UDK 519.46}
\begin{center}

{\large\bf  Separation of spectra in analysis of
Berezin kernels}

\vspace{22pt}

{\large Neretin Yu.A.\footnote{ Partially supported by
 grant RFBR 98-01-00303
and Russian program of support of scientific schools
(grant RFBR 96-01-96249).}}

 neretin@main.mccme.rssi.ru

\vspace{22pt}

\end{center}

Let $G$ be a classical real group, let  $K$ be a  maximal compact subgroup
in $G$. There exists a canonical embedding of a symmetric space
 $G/K$ to some non-compact hermitian symmetric space
  $G^\circ/K^\circ$ such that $\dim_\R G/K=\dim_\C G^\circ/K^\circ $,
%(see . tablicu v \cite{Ner2}, Dobavlenie A).
(see Addendum A). Consider a problem about  restriction of a highest
weight representation of the group   $G^\circ$
to the subgroup  $G$. We use the term "kernel-representations"
(by analogy with Bergman kernel-function) for restrictions
of this type.

A studying of kernel-representations
was initiated by F.A.Berezin
 \cite{Ber2}
shortly before his death (proofs were published only in 1994
by Upmeier and Unterberger  \cite{UU}). In last years
an interest to this subject appeared again.
Partially it is related to continuation of Berezin work,
partially this interest has independent reasons (\cite{UU},
\cite{NO}, \cite{OO}, \cite{OZ},\cite{vD},  \cite{Ner3}, \cite{MD},
\cite{Zhang}, \cite{Nerp}).
It appeared that spectra of kernel representations are very
intricate.
 For instance (see \cite{Ner4})
they include all representations which can appear in spectra
of Howe dual pairs; a decomposition of
 Howe dual pairs is known as extremely rich
spectral problem.
  G.I.Olshanskii and myself used kernel-representations
 for
construction of exotic unitary representations of the groups
 $\O(p,q),\U(p,q),
\Sp(p,q)$, see  \cite{NO}.

Emphasis that in the case  of spaces of scalar-valued functions (see \S1)
 a limit of kernel representations as a parameter
$\alpha$ tends $+\infty$ is the usual representation of the group $G$ in the space
 $L^2(G/K)$, i.e. the kernel-representations
can be considered as natural "deformation" of the spaces $L^2(G/K)$%
\footnote{In a vector-valued case, a kernel-representation can be obtained
by a deformation of a space of $L^2$-sections of a vector bundle
over $G/K$).}.

The main purpose of this paper is a decomposition of  kernel-representations
to "blocks" having more or less uniform spectra ( in fact we consider
only the case $G=\O(p,q)$; the cases  $G=\U(p,q), \Sp(p,q)$ are similar
and all information necessary for repetition of our construction for $\U(p,q)$,
$\Sp(p,q)$ is
contained in the paper). The idea to construct decompositions of this
type in various problems of harmonic analysis
is old
 (see Gelfand--Gindikin  \cite{GG}).
G.I.Olshanskii \cite{Ols1} proposed a way (noncommutative Hardy spaces)
for separation of highest weight increments
 (see survey \cite{FO} on this subject).
There are also papers of Gindikin \cite{Gin} and Molchanov \cite{Mol}
containing solution of the problem for hyperboloids.
 Our method is based on theorems about
restrictions of discontinuous functions to submanifolds.
In  \cite{Ner1}, \cite{Ols3}, \cite{NO},
\cite{vD},\cite{Ner3} this method was used for separation of the discrete part
of spectrum).

In \S1 we define kernel representations for groups $\O(p,q)$ (in the spaces
of scalar-valued functions). In \S2  we evaluate some integrals over orthogonal groups,
which seems pleasant itself. In \S3 we prove existence of restriction operators. In \S4 we prove similar statements for general kernel representations
 (in spaces of vector-valued functions).

Addendum A contains brief
discussion of definition of kernel representations

Addendum B is devoted to analytical continuation of
 Plancherel formula for kernel representations.

 I thanks G.I.Olshanskii, V.F.Molchanov, and B.Orsted   for
numerous discussion of this subject. I also thanks
H.Schlichtrull and G. van Dijk for discussions and references.

\vspace{22pt}

{\large\bf \S 1. Kernel-representations}

\vspace{22pt}

\addtocounter{sec}{1}
\setcounter{equation}{0}
\setcounter{fact}{0}

{\bf 1.0. Positive definite kernels.} Let $\cal X$ be a set. Recall that a function
$K(x,y)$  on ${\cal X}\times {\cal X}$ is called by {\it positive definite
kernel} if for all $x_1,\dots,x_n
\in{\cal X}$ the matrix
$$\left(\begin{array}{ccc}
K(x_1,x_1 )& \dots & K(x_1,x_n)   \\
\vdots&\ddots&\vdots    \\
K(x_n,x_1) &\dots &    K(x_n,x_n)
        \end{array}\right)$$
is positive definite. A positive definite kernel defines
an unique (up to a natural equivalence) Hilbert space
 $H$ and a system of vectors $\theta_x\in H$ ({\it supercomplete basis} or
{\it system of coherent states}),
enumerated by points
$x\in{\cal X}$ such that

 1. $\langle \theta_x,\theta_y\rangle_H=K(x,y)$  for all
$y,x\in{\cal X} $

2. The linear span of vectors $\theta_x$ is dense in $H$.

%Ljubomu elementu $v\in H$ my postavim v sootvetstvie funkciju $f_v$ na
%${\cal X} $ po pravilu
%$$f_v(x)=\langle v, \theta_x\rangle_H$$

For each vector $v\in H$ we define a function
$f_v(x)$ on $\cal X$ by the rule
$$f_v(x)=\langle v, \theta_x\rangle_H$$
We denote by $H^\circ$ the space of all function obtained in this way.

{\bf 1.1. Matrix balls.} Let  $p\le q$. By $\B_{p,q}$  we denote
the space of  $p\times q$  matrices with norm $<1$ over $\R$
(a norm of a matrix is the norm of operator in Euclidean space). By
$\OB_{p,q}$ we denote the closure of  $\B_{p,q}$  in $\R^{pq}$, i.e.
the set of matrices
with norm $\le 1$.

By $\O(p,q)$ we denote a group of real $(p+q)\times (p+q)$ matrices
$g={\left(\begin{array}{cc}a&b\\c&d\end{array}\right)}$
preserving bilinear form with the matrix
${\left(\begin{array}{cc}1&0\\0&-1\end{array}\right)}$.
The group $\O(p,q)$ acts on the matrix ball  $\B_{p,q}$
(and on $\OB_{p,q}$)
by fractional linear transformations
$$z\mapsto z^{[g]}:=(a+zc)^{-1}(b+zd)$$
It is easy to check that the action of $\O(p,q)$ on    $\B_{p,q}$ is transitive,
and the stabilizer of the point $z=0$ consists of matrices
${\left(\begin{array}{cc}a&0\\0&d\end{array}\right)}$
where $a\in\O(p), d\in\O(q)$.
Hence
 $$\B_{p,q} \simeq \O(p,q)/\O(p)\times\O(q)$$.

By $\M_h$ we denote the space of matrices $z\in \OB_{p,q}$
such that rank of $(1-zz^t)$ is $h$. Obviously,
the sets $\M_h$ are  orbits of the group $\O(p,q)$
 on the boundary of the matrix ball $\OB_{p,q}$. By
$\overline\M_h$ we denote the closure of the orbit   $\M_h$:
$$ \overline\M_h  =\cup_{r\le h} \M_r$$

{\bf 1.2. Berezin spaces.}
Let
$$\alpha=0,1,2,\dots,p-1 \quad\mbox{  or} \quad \alpha>p-1$$
The {\it  Berezin kernel} $K_\alpha=K_\alpha(z,u)$ is
the function on $\B_{p,q} \times \B_{p,q}$ given by the formula
$$ K_\alpha(z,u) =\det(1-zu^t)^{-\alpha}=
\det(1-uz^t)^{-\alpha}$$
By well-known Berezin theorem \cite{Ber1} for $\alpha$ indicated above
the Berezin kernel  is positive definite%
\footnote{There were also papers of Gindikin(1975), Vergne--Rossi(1976)
 and Wallach(1979) on
this subject, see references in \cite{NO} or \cite{FK}.}.
Hence we obtain the space of functions on the matrix ball.
Let us describe this space.

By ${\cal D}'(\B_{p,q})$ we denote the space of distributions
supported in the open matrix ball  $\B_{p,q}$.
Let us define the scalar product in ${\cal D}'(\B_{p,q})$
by the formula
$$\langle \phi,\psi\rangle_\alpha:= \{\phi(z)\psi(u),\det(1-zu^t)^{-\alpha}\}$$
The Berezin space $H_\alpha=H_\alpha (\B_{p,q})   $
is the completion of pre-Hilbert space
${\cal D}'(\B_{p,q})$
(in our case, a super-complete basis $\theta_x$ consists of $\delta$-functions).

{\bf 1.3. Dual realization of Berezin spaces and kernel-representations.}
By $\delta_a$ we denote the delta-function supported at the point
$a\in \B_{p,q}$.  To each element $v\in H_\alpha (\B_{p,q})$
we assign the function $f_v$ on $ \B_{p,q}$ defined by the formula
$$f_v(z) =\langle v,\delta_z\rangle_{H_\alpha}$$
Hence we obtain some space
  $H_\alpha^\circ (\B_{p,q})\simeq H_\alpha (\B_{p,q})$
consisting of functions on $\B_{p,q}$

It is readily seen that elements of the space   $H_\alpha^\circ (\B_{p,q})$
are real analytical functions on $\B_{p,q}$. Then we define the unitary representation
of the group  $\O(p,q)$ in $H_\alpha^\circ (\B_{p,q})$ by the formula
\begin{equation}T_\alpha\left(\begin{array}{cc}a&b\\c&d\end{array}\right) f(z)
=(a+zc)^{-\alpha}f\left((a+zc)^{-1}(b+zd)\right)
\end{equation}
We name the representations $T_\alpha$ by term {\it kernel-representations}

{\sc Remark.} Berezin  \cite{Ber2}  in the case of hermitian
symmetric spaces used the term {\it canonical representation}.
He gave this term  from the paper
iz \cite{VGG}. This term can have many different senses
   ; in  fact A.M.Vershik, I.M.Gelfand,
 M.I.Graev
themself in  paper \cite{VGG} and subsequent works use the term
"canonical representation" for multiplicative integrals and
and infinitely divisible representations. In fact the class of
infinitely divisible representations have very small
intersection with the class of "kernel-representations", this small
intersection is a chance.

{\bf 1.4. Holomorphic continuation.}
 By $\B_{p,q}(\C)$ we denote the set of all complex
 $p\times q$ matrices with norm $<1$.
The kernel $K_\alpha$ extends to the function
  $\det(1-zu^*)^{-\alpha}$ on
$\B_{p,q}(\C)\times \B_{p,q}(\C) $, which is holomorphic in
$z$ and antiholomorphic in $u$; this function is positive
definite on $\B_{p,q}(\C)$ (\cite{Ber1}). Then we deduce
that functions  $f\in H_\alpha^\circ (\B_{p,q}) $
extend to holomorphic functions on $\B_{p,q}(\C)$.

\vspace{22pt}

{\large\bf \S 2. Some integrals over orthogonal groups}

\vspace{22pt}

\addtocounter{sec}{1}
\setcounter{equation}{0}
\setcounter{fact}{0}

Consider a matrix $A$. By $[A]_p$
we denote its left  upper corner having a size $p\times p$.
   We denote the Haar measure on $\SO(n)$ by
 $d\sigma_n(g)$.

{\bf 2.1. The maps $\Upsilon^m$.}
 Let us wright a matrix $g\in\SO(n)$ as a $(m+(n-m))\times (m+(n-m))$
block matrix
 $\left(\begin{array}{cc}P&Q\\R&T\end{array}\right)$.
 Consider the map
$$\Upsilon^m:
\left(\begin{array}{cc}P&Q\\R&T\end{array}\right) \mapsto T-R(1+P)^{-1}Q$$
defined almost everywhere

{\sc Proposition 2.1.}{\it
  a) $\Upsilon^m$ is a map $\SO(n)$ to $\SO(n-m)$.
\begin{eqnarray*}
 \mbox{ b})& \Upsilon^k\circ\Upsilon^m=\Upsilon^{k+m}.\\
 \mbox{ c})& S=\frac{g-1}{g+1}\,\mbox{ implies }
\{S\}_p=\frac{\Upsilon^{n-p}(g)-1}{\Upsilon^{n-p}(g)+1}
\end{eqnarray*}
where $\{S\}_p$ denotes the  right lower $p\times p$ corner of the matrix  $S$.
 }

{\sc Remark.} The map $\Upsilon$ is a value of Livshic
characteristic function
$\chi(z)=
 T+zR(1-zP)^{-1}Q$
at the point $z=-1$ (see \cite{Liv}).

{\sc Proof.} Statement c) can be easily checked by
 direct calculation using Frobenius formula for inversion
of block matrix  (see \cite{Gan}, \S {
\rm II}.4)
,
$$\left(\begin{array}{cc}A&B\\C&D\end{array}\right)^{-1}=
\left(\begin{array}{cc}
A^{-1}+A^{-1}B(D-CA^{-1}B)^{-1}C A^{-1}&
- A^{-1}B(D-CA^{-1}B)^{-1}
\\-(D-CA^{-1}B)^{-1} C A^{-1} & (D-CA^{-1}B)^{-1}
\end{array}\right)$$

Statements a),b)   are consequences of c).

{\sc Proposition 2.2.} {\it Let $A,B,\in \SO(n-m)$. Then
$$\Upsilon^m\left[\left(\begin{array}{cc}1&0\\0&A\end{array}\right)   g
\left(\begin{array}{cc}1&0\\0&B\end{array}\right)\right]=
A\Upsilon^m(g)B$$               }

{\sc Proof}  is obvious.

{\bf 2.2. Another proof of Proposition 2.1.}
Consider the space $V_n:=\R^n\oplus\R^n$ equipped with indefinite
bilinear form
$$L_n((x,y),(x',y')):=\sum_{k=1}^n x_k x'_k -
\sum_{k=1}^n y_k y'_k$$
A graph of an unitary operator  $\R^n\to\R^n$
is a $L_n$-isotropic subspace in
 $V_n=\R^n\oplus\R^n$.

For each linear relation   $M:V_n\rightrightarrows V_{n-m}$
(recall that a linear relation is
a subspace $M\subset V_n\oplus V_{n-m}$, for details see \cite{Ner2})
and each subspace $H\subset V_n$   it is possible to define
a subspace $MH\subset V_{n-m}$ in the following way. A vector
$y\in V_{n-m}$ is an element of $MH$ iff there exists  $x\in H$ such that
$(x,y)\in M$.

Consider the linear relation
 ${\cal Z}_n^{n-m}:V_n\rightrightarrows V_{n-m}$
consisting of all vectors having the form
$$\bigr(\bigl\{ (v,h),(w,-h)\bigr\},\bigl\{v,w\bigr\}\bigr)\in
\bigl\{\R^n\oplus\R^n\bigr\}\oplus \bigl\{\R^{n-m}\oplus\R^{n-m}\bigr\}$$
The condition $(u,v),(u',v')\in {\cal Z}_n^{n-m}$ implies $L_n(u,u')=L_{n-m}(v,v')$.
Hence the transformation   $H\mapsto {\cal Z}_n^{n-m}H$
maps isotropic subspaces to isotropic subspaces.

It is readily seen that the map $\Upsilon^m$ coincides with
the transformation
 $H\mapsto {\cal Z}_n^{n-m}H$ of isotropic grassmannians.
Now Proposition 2.1 becomes obvious.

{\bf 2.3. A projection of the Haar measure.}
 Denote by $I$ the segment $[-1,1]$.
Consider the map
 $$\Xi:\SO(n)\to  \SO(n-1)\times I$$
given by the formula (see a notation in the beginning of the section)
$$\Xi(g)=\left( \Upsilon^1(g), [g]_1\right)$$

{\sc Proposition 2.3.} {\it The image of the Haar measure $d\sigma_n$
 under the map $\Xi$ is
$$\const\cdot(1-x^2)^{(n-3)/2}dx\,d\sigma_{n-1}$$}

{\sc Proof.} By $\SO(n-1)$- equivariance of the map
$\Upsilon^1$ (Proposition 2.2), the image of the Haar measure
has the form
$\phi(x) dx\,d\sigma_{n-1} $. On another side the projection of the Haar measure
to $I$ is $\const\cdot (1-x^2)^{(n-3)/2}$ and we obtain required form
of the function
$\phi(x)$.

{\bf 2.4. The map of the orthogonal group to the cub.}
Consider the maps
$$\SO(n)\to  \SO(n-1)\times  [-1,1]
\to  \SO(n-2)\times[-1,1] \times  [-1,1]\to\dots$$
 As result we obtain the map of the group $\SO(n) $ to the cube
$[-1,1]^{n-1}$,
given by the formula
$$(x_1,x_2,\dots,x_{n-1})=\bigl
([\Upsilon^{n-2}(g)]_1, \dots,
[\Upsilon^1(g)]_1,[g]_1\bigr)$$

%Otmetim, chto gruppa $\O(1)$ sostoit iz dvuh tochek $+1,-1$,
%poetomu kazhdaja iz komponent svjaznosti gruppy $\O(n)$ pri etom
%otobrazhaetsja na kub  $[-1,1]^{n-1}$.

By Proposition 2.3, the image of the Haar measure under this map is
$$d\mu(x_1,x_2,\dots,x_{n-1})=\const\cdot\prod_{j=1}^{n-1} (1-x_j^2)^{(j-2)/2}
d\nu(x_1)\prod_{j=1}^{n-1}dx_j$$
%where $\nu$ -- mera na dvoetochii $x_1=\pm1$ takaja, chto mera kazhdoj tochki
%ravna $1/2$.
Hence for each function $f$ depending on $n-1$ variables
 $x_1,\dots, x_{n-1}\in [-1,1]$
 we have
\begin{eqnarray}
\int_{\SO(n)}
f([\Upsilon^{n-1}(g)]_1, \dots,[\Upsilon^1(g)]_1,[g]_1)\,d\sigma_n(g)=
\nonumber \\
=\int_{[-1,1]^{n-1}} f(x_1,x_2,\dots,x_{n-1})  d\mu(x_1,x_2,\dots,x_{n-1})
\end{eqnarray}

{\bf 2.5. Multiplicativity.}
{\sc Proposition 2.4.}{\it Let $g\in\O(n)$. Let $m<p\le n$. Then}
$$\det(1+[g]_p)=\det(1+[g]_m)\det(1+[\Upsilon^m(g)]_{p-m})$$

{\sc Proof.} Let us wright $g$  as block
$(m+(p-m)+(n-p)) \times (m+(p-m)+(n-p)) $ matrix:
$$g=\left(\begin{array}{ccc}P&Q_1&Q_2\\R_1&T_{11}&T_{12} \\
R_2&T_{21}&T_{22}\end{array}\right) $$

By the usual formula for determinant of a block matrix, we obtain %(\cite{Gan})
$$\det\left[1+\left(\begin{array}{cc}
P&Q_1\\R_1&T_{11}\end{array}\right)\right]=
\det(1+P)\det(1+T_{11}-R_1(1+P)^{-1}Q_1)$$
On another side
$$\Upsilon^m(g)=
\left(\begin{array}{cc}T_{11}-R_1(1+P)^{-1}Q_1&T_{12}-R_1(1+P)^{-1}Q_2  \\
      T_{21}-R_2(1+P)^{-1}Q_1  & T_{22}-R_2(1+P)^{-1}Q_2\end{array}\right) $$
and now the statement becomes obvious.

{\bf 2.6. Some integrals.}
{\sc Theorem 2.5.} {\it Assume that the Haar measure of the whole group $\SO(n)$
is 1. Let $\lambda_1,\dots,\lambda_n\in\C, \lambda_{n+1}=0$, and
${\rm Re}\,\lambda_k>-(n-k)/2$
for all  $k=0,1,\dots,n$. Then
\begin{eqnarray}\int_{\SO(n)}\prod_{k=1}^{n-1}
\det(1+[g]_{k})^{\lambda_k-\lambda_{k+1}}\,d\sigma_n(g)=
% \nonumber\\
\prod_{k=1}^{n-1}\frac{\Gamma(n-k)\Gamma(\lambda_k+(n-k)/2)}
{\Gamma((n-k)/2)\Gamma(\lambda_k+n-k)}
\end{eqnarray}      }

{\sc Proof.} Let us apply  formula (2.1) and Proposition 2.4. Then our
 integral splits to the product of integrals having the form
\begin{eqnarray}
& \int_{-1}^1(1-x^2)^{(n-k-2)/2}(1-x)^{\lambda_k} dx=
\\&= 2^{{\lambda_k}+n-3}B\left({\lambda_k}+\frac{n-k}{2},\frac{n-k}{2}  \right)=
 2^{{\lambda_k}+n-k-3} \frac{\Gamma({\lambda_k}+(n-k)/2) \Gamma((n-k)/2)}
{\Gamma({\lambda_k}+n-k)}\nonumber
\end{eqnarray}

{\bf 2.7. Some integrals over unitary and symplectic groups.}
The following integrals can be evaluated in the same way, the only
the obvious integral
(2.3) is replaced by more complicated integrals over ball $|x|<1$
in $\C$ or in $\H$.
\begin{eqnarray}
\int_{\U(n)}\prod_{k=1}^{n}
\det(1+[g]_{k})^{\lambda_k-\lambda_{k+1}}
\overline{\det(1+[g]_{k})}^{\mu_k-\mu_{k+1}}\,d\sigma_n(g)
=\nonumber\\=
\prod_{k=1}^n \frac{\Gamma(n-k+1)\Gamma(n-k+1+\lambda_k+\mu_k)}
{\Gamma(n-k+1+\lambda_k)\Gamma(n-k+1+\mu_k)}\\
\int_{\Sp(n)}\prod_{k=1}^{n}|\det(1+[g]_{k})|^{\lambda_k-\lambda_{k+1}}\,
d\sigma_n(g)
=\nonumber\\=
 \prod_{k=1}^n
\frac{\Gamma(2(n-k+1))\Gamma(2(n-k+1)+\lambda_k+1)}{\Gamma(2(n-k+1)+\lambda_k/2)
\Gamma(2(n-k+1)+\lambda_k/2+1)}
\end{eqnarray}
(we define the compact symplectic group  $\Sp(n)$ as the group of
unitary operators
in $n$-dimensional quaternionic space ${\Bbb H}^n$;
a quaternionic operator $A$ can be considered as an operator $A_\R:\R^{4n}\to \R^{4n}$; we define a determinant of the quaternionic operator $A$
as $\sqrt[4]{\det A_\R})$.

{\bf 2.8. Other ways of evaluation of integrals (2.2), (2.4), (2.5).}
We obtain a verbal evaluation of these integrals, nevertheless
  our way is not
quite usual. There exist another way of evaluation of the same integrals.
The Cayley transform reduces (2.2) to the form
\begin{equation}
\int \prod_{k=2}^n \det (1-[T]_k)^{\mu_k}\,dT
\end{equation}
where integration is given over the space of all skew-symmetric matrices.
Last integral can be easily evaluated by the method described in
\cite{Nerp} (in fact (2.6) is simpler than integrals from \cite{Nerp}).
Similar way is valid for (2.4),(2.5).

Integral (2.4) reduces to the form
$$\int\prod_{k=1}^n \det(1+[T]_k)^{\mu_k} \det(1-[T]_k)^{\nu_k}\,dT$$
where integration is given over the space of all anti-hermitian matrices
($T=-T^*$). In this case where are two additional possibilities: we can use
Laplace transform or we can use Gindikin method of separation of variables
(see \cite{FK}, chapter 7).

{\bf Remark.} If $\lambda_1=\lambda_2=\dots=\lambda_n$, then integral   (2.2)
is equivalent to one of Hua Loo Keng integrals \cite{Hua}, chapter 2

\vspace{22pt}

{\large\bf \S 3. Operators of restriction to the boundary}

\vspace{22pt}

\addtocounter{sec}{1}
\setcounter{equation}{0}
\setcounter{fact}{0}

In this section we assume $p< q$.

{\bf 3.1. Operators of restriction to $\M_r$.} Denote by
 $\H^\circ_\alpha$
the subspace
in $H^\circ_\alpha (\B_{p,q})$ consisting of functions analytical in
some neighborhood of closed matrix ball
 $\OB_{p,q}$.

{\sc Lemma 3.1.}{\it  The subspace  $\H^\circ_\alpha$ is dense
in  $H^\circ_\alpha $.}

{\sc Proof.} There exists a basis in $H^\circ_\alpha $
consisting of orthogonal homogeneous polynomials (see for instance\cite{NO},
\S1). These functions are eigenfunctions of the operator $A_{\epsilon}f(z)=f((1-\epsilon)z)$.
The eigenvalues are $(1-\epsilon)^k$. Hence for each $\epsilon>0$
the operator $A_\epsilon$ is a continuous   operator
$H^\circ_\alpha \to H^\circ_\alpha$. Obviously the image of the operator
$A_\epsilon$
is dense in  $H^\circ_\alpha $ and also it is contained in
  $\H^\circ_\alpha$.
 $\blacksquare$

By $\L(\M_r)$ we denote the space of locally integrable functions
on manifold $\M_r$ (see Subsection 1.1;
if $r\ne0$ then $M_r$ is not compact).

The group $\O(p,q)$ acts in the space $\L(\M_r)$ by same formula (1.1),
as in $H^\circ_\alpha(\B_{p,q})$
(we only have to think that  a matrix
$z$ is an element of $\M_r$).

By  $J_r$ we denote the operator of restriction of
a function $f\in \H^\circ_\alpha$
to the submanifold $\M_r$.

{\sc Theorem 3.2.} {\it  Let $\alpha<(q-p+2r)/2$. Then the operator $J_r$
extends to a continuous operator $H^\circ_\alpha\to \L(\M_r) $.}

\smallskip

{\sc Remark.} Let $\alpha\ne 0$ (if $\alpha=0$ then the space
$H^\circ_\alpha$ is one-dimensional).  In general, functions $f\in H^\circ_\alpha$,
are discontinuous in point of boundary of the matrix ball. In other words,
for points $a\in \OB_{p,q}\setminus \B_{p,q}$
the lineal functional $f\mapsto f(a)$ on $H^\circ_\alpha$
is not well-defined.
Hence existence of the restriction operator is not
obvious.

\smallskip

Evidently, the operator $J_r: H^\circ_\alpha\to \L(\M_r) $ is
 $\O(p,q)$-intertwining. Hence its kernel $\ker J_r$ is a
subrepresentation.
Hence the orthogonal complement  $(\ker J_r)^\bot$ also is a subrepresentation.
The operator $J_r$ is a bijection from $(\ker J_r)^\bot$
to the image  ${\rm im}\, J_r$ of the operator $J_r$.
Hence we constructed a subrepresentation in
$H^\circ_\alpha$,
and this subrepresentation has  a natural realization in a space
of functions on $\M_r$.

{\bf 3.2. Proof of the theorem 3.2.} In \cite{NO} (theorems 2.2--2.2)
 there was obtained
one general statement on boundary values of holomorphic functions.
In our case, it implies the following statement

{\bf Theorem 3.3.} {\it Let $M$ be a compact
subset in the boundary of the matrix ball $\B_{p,q}$.
Let $\mu$ be a measure  supported in  $M$.
 Assume that

i) For almost all (with respect to the measure $\mu\times\mu$)
points $(z,u)\in M\times M$
there exists a limit
as $\epsilon$ tends to 0 from the right
$$L(z,u)=\lim_{\epsilon\to+0}\det(1-(1-\epsilon)zu^t)^{-\alpha}$$

ii) This limit is dominated, i.e. there exists a function
 $h\in L^1(\mu\times\mu)$ such that
$$\det(1-(1-\epsilon)zu^t)^{-\alpha} \le h(z,u)$$

Then the restriction operator of a function $f\in H^\circ_\alpha$ to the subset $M$
is a continuous  operator  $H_\alpha^\circ \to L^1(M)$.}

 First we will explain briefly a way of application of the theorem
in our case.
We want to apply it to the orbit $\M_r$. Obviously the limit exists.
In fact it is sufficient to check the convergence of integral
of the function $L(z,u)=\det(1-zu^*)^{-\alpha}$.
 Nevertheless I don't see a natural measure on
  $\M_r$.
Hence there is a small hope to evaluate integral if we don't know
natural integrand.

An estimation of convergence also seems heavy since the singularity
of the integral is very unpleasant. In fact we construct
 non-directly a calculable integral with given singularities.
w

\medskip

To apply Theorem 3.3,
we consider arbitrary compact subset  $M$ on
the manifold $\M_r$,
we also assume that $M$ is a closure of its interior.
Let $\mu$ be the surface Lebesgue measure
 on $M$. It is sufficient to check that the integral
\begin{equation}\int_M\int_M \det(1-zu^t)^{-\alpha}d\mu(z)\,d\mu(u)\end{equation}
is convergent.  The condition  i) is a consequence of the convergence.
Condition ii) follows from the obvious inequality
(\cite{NO},3.8):
$$\det(1-c^2 zu^t)^{-\alpha}\le 2^{p\alpha}\det(1-zu^t)^{-\alpha}.$$

It is sufficient to check that the integral
\begin{equation}\int_M \det(1-zu^t)^{-\alpha}d\mu(z)\end{equation}
is convergent for any $u$ and is dominated by a constant
that don't depend on $u$. As it will be shown below, it is sufficient
to prove convergence of the integral for some  $u$.

First we consider $u=u^{(0)}$ given by
block $((p-r)+r)\times((p-r)+(q-p+r))$ matrix
$\left(\begin{array}{cc}-1&0\\0&0\end{array}\right)$.
Then  integral (3.2) transforms to the form
\begin{equation}\int_M \det[1+z]_{p-r}^{-\alpha}d\mu(z)\end{equation}

{\sc Lemma 3.4.} {\it Integral} (3.3) {\it converges.}

{\sc Lemma 3.5.} {\it Let us represent $g\in\SO(q+r)$ as block matrix
$g=
\left(\begin{array}{cc}g_{11}&g_{12}\\g_{21}&g_{22}\end{array}\right)$
having the size $(p+(q+r-p))\times(q+r)$. Then $g_{11}\in \OM_r$. Conversely each
matrices $z\in\OM_r$ is a left upper corner of some matrix $g\in\SO(q+r)$.}

{\sc Proof of Lemma 3.5} is obvious.

{\sc Proof of Lemma 3.4}. By Lemma 3.5 the map
 $\tau:g\mapsto g_{11} $ takes  $\SO(q+r)$
to $\OM_r$. Obviously, the image of the Haar measure on $\SO(q+r)$
under the map $\tau$ is a measure equivalent to
a surface Lebesgue measure on $\M_r$.
Hence it is sufficient to check convergence of the integral
$$\int_{\SO(q+r)}\det[1+g]_{p-r}^{-\alpha}d\sigma_{q+r}(g)$$
But this integral was evaluated in \S2
( we substitute $\lambda_1=\dots=\lambda_{p-r-1}=-\alpha; \lambda_{p-r}=\dots=0$ to
formula (2.2)).                                       $\blacksquare$

{\sc Lemma 3.6.} {\it
 The kernel
$\det(1-zu^t)^{-\alpha}$   satisfies to the equality         }
\begin{equation}\det(1-z^{[g]}(u^{[g]})^t)^{-\alpha}=
\det(1-zu^t)^{-\alpha} \det(a+zc)^\alpha\det(a+zu)^\alpha
\end{equation}

{\sc Proof} is obvious.

Consider a compact subset  $L\subset\O(p,q)$
such that for each $v\in\M$ there exists $g_v\in L$ satisfying the condition
$v^{[g_v]}=u^{(0)}$. Denote by $K$ the set of all points
in $\M_r$ having the form $u^{[g]}$ where $u\in M$ and $g\in L$.
Obviously, $K$ is compact.

Substitute $w=z^{[g]}$ to the integral (3.2). By formula (3.4),
we obtain
\begin{equation}
\int\limits_{M^{[g_u]}} \det(1-u^{(0)}w^t)^{-\alpha} H(w,u)\,d\mu(w),
\end{equation}
where $M^{[g_u]} $ denote the image of the set $M$ under the map $g_u$,
and $H(w,u)$ is some function (a product of Jacobian and
additional factor which appears from (3.4)).
Obviously, the function
$H(w,u)$ is bounded by some constant  $C$.
Hence integral (3.5) is dominated by quantity
$$C\int_K  \det(1-u^{[0]}w^t)^{-\alpha} d\mu(w)  $$
This quantity doesn't depend on  $u$ and it is finite by  Lemma 3.4.
Theorem is proved.

{\bf 3.3. Restriction of partial derivatives to the boundary.}
In Subsection 3.2 we restricted functions to $\M_r$.
Also it is possible to restrict partial derivatives of the function $f\in H^\circ_\alpha  $. The existence of restriction operator for
partial derivative $\frac{\partial}{\partial z_{ij}}$
depends on convergence of integral
$$\int_{M} \int_{M} \left|
\frac{\partial}{\partial z_{ij}}  \frac{\partial}{\partial u_{ij}}
\det(1-zu^t)^{-\alpha}\right|\,dz\,du,$$
where $M\subset \M_r$ is a compact subset.
The integrand has the form
$|P(z,u)|\det(1-zu^t)^{-\alpha-2}$ where $P(z,u)$ is a polynomial.
Hence we can apply the same arguments.
As a result, we obtain the following theorem.

{\sc Theorem 3.7.} {\it Let $\alpha+2k<(q-p+2r)/2$.
Let $\partial^\nu$
be a partial derivative of order  $k$. Then the operator of restriction
of a function  $\partial^\nu f(z)$ to $\M_r$ is a continuous
operator
from $H_\alpha^\circ$ to $\L(\M_r)$.
}

Now we will slightly formalize the picture.

Denote by $\NN(\M_r)$ the conormal bundle to $\M_r$ in $\R^{pq}$.
By  $\SS^k\NN(\M_r)$ we denote  $k$-th symmetric power
of the conormal bundle. By $\F_r^k$ we denote the set of all
 $f\in \H^\circ_\alpha$ whose partial derivatives
of order $\le k$ vanish on $\M_r$.
Then elements of quotient space
$\F_r^{k-1}/\F_r^{k}$ are identified with sections of the
bundle $\SS^k\NN(\M_r)$.

Theorem 3.7 implies the following statement:

 {\it Let $\alpha+2k<(q-p+2r)/2$.
Then the operator of restriction of  partial derivatives of order $k$ to
$\M_r$ is a well-defined operator from the space
 $H_\alpha^\circ$ to the space
of $L^1_{\rm loc}$-sections of the bundle
  $\SS^k\NN(\M_r)$.}\footnote{ the algebraic version of this
construction arises to \cite{JV}}

{\bf 3.4. Some comments.}  Thus we obtain the intertwining
operator
$J_{r,k}$
from  $H_\alpha^\circ$ to the space of sections
of bundle $\SS^k\NN(\M_r)$.
Its image  ${\rm im} J_{r,k}$  is
a quotient space of a Hilbert space.
Hence ${\rm im} J_{r,k}$ has a natural structure of a Hilbert space.
(we emphasis that the Hilbert topology in
  ${\rm im} J_{r,k}$ is stronger than topology induced from
$\L$). We denote the representation of the group $\O(p,q)$ in
 ${\rm im} J_{r,k}$
by
\begin{equation}\Xi_{r,k}^\alpha=\Xi_{r,k}^\alpha(p,q)\end{equation}
The quotient space of a Hilbert space  $H$ by arbitrary
closed subspace $L$ can be canonically identified with
the orthogonal complement
$L^\bot$ to $L$. Hence the representation  $\Xi_{r,k}^\alpha$ is
a subrepresentation
in the kernel-representation $T_\alpha$.

Hence we obtained a canonical family of subspaces (3.6) in
 $H^\circ_\alpha(\B_{p,q})$, we name these subspaces by the
term {\it blocks}.

\smallskip

If $\alpha$ is sufficiently large (i.e. $\alpha >(p+q)/2-1 $), then
the operators of restriction to submanifolds $\M_r$ are not defined.
In this case blocks  $\Xi_{r,k}^\alpha(p,q)$ don't exist.
It is known that for
 $\alpha >(p+q)/2-1 $ the representation of $\O(p,q)$
in $H^\circ_\alpha(\B_{p,q})$ is equivalent to the representation of  $\O(p,q)$ v
in $L^2$ on symmetric space  $\B_{p,q}$\footnote{this statement
for arbitrary scalar-valued kernel-representations is proved
in cite  \cite{OO}; earlier similar statement for kernel representations
of some groups were obtained (in different degree of generality) in
\cite{Ber2}, \cite{Rep}, \cite{Gut}).}.

\smallskip

For all $\alpha>p-1$ the summand ("block"), that equivalent to $L^2(B_{p,q})$,
is present in decomposition of the kernel-representation.

For $\alpha <(p+q)/2-2 $
the block $\Xi^\alpha_{p-1,0}$ appears. This summand corresponds
to the operator of restriction of function to the maximal boundary orbit $M_{p-1}$.

For $\alpha<(p+q)/2-3$ block $\Xi^\alpha_{p-2,0}$ appears.
After point $\alpha=(p+q)/2-4$ two additional blocks $\Xi^\alpha_{p-1,1}$
and $\Xi^\alpha_{p-3,0}$ appears etc.

\smallskip

At the point  $\alpha=p-1$ (it is the end of continuous series, see 1.2)
the largest block $L^2 (\B_{p,q})$ disappears from decomposition.
Also at this point all blocks $\Xi_{p-1,1}^{p-1}$, $\Xi_{p-1,2}^{p-1}$, $\Xi_{p-1,3}^{p-1}$ etc. disappear. But the block $\Xi_{p-1,0}^{p-1}$
at this point still alive.

If we continue decreasing of $\alpha$ (after point $\alpha=p-1$ it
is a discrete process, see Subsection 1.2)
"new" blocks continue their appearance. Nevertheless
"old" blocks quickly disappear (this process can be checked using
Plancherel formula, see Addendum B).

\smallskip

 At the point
 $\alpha=0$ only the block  $\Xi^0_{0,0}$ is still alive, and this
block is one-dimensional representation.

\smallskip

Let us briefly discuss spectra of representations  $\Xi^\alpha_{r,k}$.

\smallskip

%Oboznachim cherez $P_r$ parabolicheskuju podgruppu v $\O(p,q)$, javljajuschujusja stabilizatorom
%$(p-r)$-mernogo izotropnogo podprostranstva. Ocheidno reduktivnaja chast' gruppy
%$P_r$ est' $\GL(p-r,\R)\times\O(r,q-p+r)$.

{\sc Proposition 3.8.} {\it Each subrepresentation in the kernel-representation
 $T_\alpha$ contains an $\O(p)\times\O(q)$-invariant vector.}

{\sc Proof.} The function $f(z)=1$ is an $\O(p,q)$-cyclic
 vector in $H^\circ_\alpha$
(since its orbit is exactly super-complete basis).
 Hence the projection of the function    $f(z)=1$
to any invariant subspace is a cyclic vector in this subspace.
Hence the projection is non zero.

{\sc Proposition 3.9.}
{\it {\rm a)} The representation of the group $\O(p,q)$ in
 $\Xi^\alpha_{0,0}$
is irreducible

{\rm b)} The representation of  $\O(p,q)$ in $\Xi^\alpha_{0,k}$
is a finite sum of irreducible representations.
}

{\sc Proof.} a) Indeed  the constant is  the unique
$\O(p)\times\O(q)$-invariant function on $M_0$.

b) We have to check that the space of
$\O(p)\times\O(q)$-invariant sections of
the bundle  $\SS^k\NN(\M_0)$ is finite-dimensional.
But dimension of this space can't exceed the dimension
of fiber.

{\sc Remark.} "Block" $\Xi^\alpha_{0,0}$ which corresponds to
operator of restriction to the compact orbit $\M_0$
was discovered in the paper of G.I.Olshanskii and myself
\cite{NO}, \S 3
(under more rigid conditions to parameters of representations) $\alpha<(q-2p+1)/2$.

{\sc Remark.} Spectrum in blocks  $\Xi^\alpha_{r,k}$ (where $0<r<p$)
is purely continuous. It consist of representations
induced from unitary representations of groups
$\O(p-r,q-r)$. The last representations are contained in
spectrum of
$\Xi^{\alpha+2k}_{0,k}(p-r,q-r)$.
%This statement is a consequence
%(and in a parallel way an element of the proof) of the Plancherel formula
%(see Addendum B).

{\sc Remark.} "Additional"   block with continuous spectrum in kernel-representations
was discovered earlier only for the case  $G=\U(2,2)$  in the paper \cite{OZ}.
In \cite{Ner4} I conjectured that the spectrum of kernel-representations
for small values of $\alpha$ contains many various series.
The proof of this conjecture is the subject of this paper.

\vspace{22pt}

{\bf \S 4. General kernel-representations.}

\vspace{22pt}

{\bf 4.1. Definition.}  By $\rho_\Lambda$ we denote finite dimensional
holomorphic representation
 of the group   $\GL(n,\C)$ defined by a signature
$\Lambda=(\lambda_1, \dots,\lambda_n)$, where
 $\lambda_1\ge \dots\ge\lambda_n$. By $V_\Lambda$ we denote
the representation $\rho_\Lambda$.
Recall that there exists a canonical $\U(n)$-invariant scalar product
in $V_\Lambda$.

Consider the representations $\rho_{\Lambda'}$, $\rho_{\Lambda''}$
of the groups  $\GL(p,\C)$, $\GL(q,\C)$ defined by signatures
$\Lambda'=(\lambda'_1, \dots,\lambda'_p)$,
$\Lambda''=(\lambda''_1, \dots,\lambda''_q)$ respectively.
We will use notations
 $\widetilde\Lambda= (\lambda'_1, \dots,\lambda'_p;
\lambda''_1, \dots,\lambda''_q)$ and also
$V_{\widetilde\Lambda}=V_{\Lambda'} \otimes
V_{\Lambda''}            $.

Consider the function
$$ K_{\widetilde\Lambda}(z,u)=
\rho_{\Lambda'}(1-zu^t)\otimes \rho_{\Lambda''}(1-uz^t);
\qquad z,u\in \B_{p,q}$$
taking values in the space  ${\rm End}(V_{\widetilde\Lambda})$
 of linear operators in
 $V_{\widetilde\Lambda}$.

Assume that $ K_{\widetilde\Lambda}(z,u)$
is a matrix-valued positive defined kernel.
Recall that this means positive definiteness
of the kernel
$$L\left((z,\xi);(u,\eta)\right)=
\langle K_{\widetilde\Lambda}(z,u)\xi,\eta\rangle_{V_{\widetilde\Lambda}};
\qquad x,y\in \B_{p,q};\quad \xi,\eta\in V_{\widetilde\Lambda}$$
on the space ${\cal X}=\B_{p,q}\times  V_{\widetilde\Lambda}$
(for the condition of positive definiteness of the kernel ${\widetilde\Lambda}(z,u)$
we obtained in \cite{Ols00}, see also any text on classification of highest weight unitary representations, see also \cite{NO}.
 %sm., naprimer, v \cite{NO}).

The kernel $L\left((z,\xi);(u,\eta)\right)$ defines a Hilbert space
consisting of functions on the space
${\cal X}=\B_{p,q}\times  V_{\widetilde\Lambda}$.
But these functions are linear in the variable  $\xi\in V_{\widetilde\Lambda}$.
Hence we can consider them as $V_{\widetilde\Lambda}$-valued
functions on $\B_{p,q}$.

As a result
 we obtain the space $H^\circ_{\widetilde\Lambda}$,
which consists of $V_{\widetilde\Lambda}$-valued functions on $\B_{p,q}$.
The group  $\O(p,q)$ acts on $H^\circ_{\widetilde\Lambda} $ by
unitary operators
given by the formula
$$T_{\widetilde\Lambda}(g) f(z)=
\left[\rho_{\Lambda'}(a+zc)\otimes\rho_{\Lambda''}(-d+cz^{[g]})\right]
f(z^{[g]})$$

{\bf 4.2. Operators of restriction to the boundary.}
{\sc Theorem 4.1.}{\it a) Assume $\lambda''_1-\lambda'_p <(q-p+2r)/2$.
Then the operator of restriction to an orbit $M_r$ is a well-defined
operator in   $H^\circ_{\widetilde\Lambda}$.

 b) Assume $2k+\lambda''_1-\lambda^{'}_p <(q-p+2r)/2$.
Then the operator of restriction of partial derivatives of order
$k$ to the orbit $M_r$ is a well-defined operator in   $H^\circ_{\widetilde\Lambda}$.
}

{\sc Proof} is similar to the proof of Theorem  3.2,
it is necessary also use simple considerations from \cite{NO},4.5.

{\bf 4.3. Comments.} As above, we obtain
 decomposition of kernel-representation to "blocks" $\Xi_{r,k}^\Lambda$,
 where  $r$ is a parameter enumerating orbits, and  $k$ is an order
of partial derivatives.

Existence of  operator of restriction to  $\M_0$ was discovered in
 \cite{NO},\S 4 (under more rigid restriction to the signature $\Lambda$).
 In particular this construction gives a simple way
for obtaining of "exotic" unitary representations of
the groups $\O(p,q)$ (see \cite{NO},\S\S 4-5).
Emphasis also that this representations admit a limit as
  $q\to\infty$ (\cite{NO}, \S 6).

\smallskip

 The representations of $\O(p,q)$ in   $\Xi_{0,0}^\Lambda$ are reducible;
their decomposition is discussed in  \cite{NO} , \S 5.
As above,
the representation of $\O(p,q)$ in all spaces  $\Xi_{0,k}^\Lambda$
has a purely discrete spectrum.
Nevertheless the spectrum of blocks  $\Xi_{r,k}^\Lambda$
for $r>0$ is  still sophisticated and  contains representations
of different types.
Here I don't have likely conjectures.

 \begin{center}

{\large\bf Addendum A. Hermitization of symmetric spaces.}

\renewcommand{\theequation}{A.\arabic{equation}}
\addtocounter{fact}{7}
\setcounter{equation}{0}

\end{center}

\renewcommand\H{{\Bbb H}}

 Consider the following list.

\def\table{\begin{tabbing}
10.\=\qquad\qquad xxxxxxxxxxxxxxxxxxxxxxxxx \= \kill}
\def\etable{\end{tabbing}}
\def\GL{{\rm GL}}
\def\G{{$G^\circ$}}
\def\Sp{{\rm Sp}}
\def\SOS{{\rm SO^*}}

\table
%\begin{tabbing}xxxxxx\= ZZZZZZZZZZZZ\=zzzzz \kill
1.\>$G=\GL(n,\R)$ \>  $G^\circ=  \Sp(2n,\R)$\\
2.\>$G=\O(p,q) $    \>  $G^\circ= \U(p,q)$\\
3.\>$G=\Sp(2n,\R)$ \>  $G^\circ= \Sp(2n,\R)\times\Sp(2n,\R)$\\
4.\>$G=\GL(n,\C)$ \> $ G^\circ=\U(n,n)$\\
5.\>$G=\SO(n,\C)$ \> $G^\circ= \SOS(2n)$\\
6.\>$G=\Sp(2n,\C)$ \> $G^\circ= \Sp(4n,\R)$\\
7.\>$G=\U(p,q) $ \> $G^\circ=\U(p,q)\times\U(p,q)$\\
8.\>$G=\GL(n,\H)$ \> $G^\circ=\SOS (2n)$\\
9.\>$G=\Sp(p,q)$ \> $G^\circ=\U(2p,2q)$\\
10.\>$G=\SOS(2n)$\>$G^\circ=\SOS(2n)\times\SOS(2n)$\\
11.\>$G=\SO(2,n)$\>$G^\circ=\SO(2,n)\times\SO(2,n)$\\
12.\>$G=\SO(1,p)\times\SO(1,q)$ \> $G^\circ=\SO(2,p+q)$
\etable

Denote by $K$, $K^\circ$
the maximal compact subgroup in $G$, $G^\circ$ respectively.
Then $G^\circ/K^\circ$ is a hermitian symmetric space and
$G/K\subset G^\circ/K^\circ $ is a totally real subspace,
and
$$\dim_\R G/K=\dim_\C G^\circ/K^\circ$$
We say that $G^\circ/K^\circ$ is {\it hermitization} of $G/K$.
If $G/K$ is hermitian then its hermitization is $G/K\times G/K$.
We see that all classical riemann symmetric spaces have hermitization%
\footnote{Unfortunately I don't know the author of this observation.
First time I heard this $20$ years ago.}

\smallskip

We define kernel-representation as restriction of an unitary
highest weight representation\footnote{A highest weight representation
of the group $G^\circ$ is a representation in a space of holomorphic vector-valued
functions on
  $G^\circ/K^\circ$}
$\rho$ of the group $G^\circ$ to the subgroup $G$%
\footnote{references:
\cite{Ber2} (the cases $G=\U(p,q), \Sp(2n,\R), \SOS(2n),\SO(2,n)$),
\cite{NO} ($G=\O(p,q),\U(p,q), \Sp(p,q)$),
 \cite{OO}, \cite{Ner4} (some discussion of other cases);
there exists an idea to extend a definition of kernel representation
to pseudoriemann symmetric spaces, see \cite{MD}.}

          \pagebreak

\renewcommand{\theequation}{A.\arabic{equation}}
\addtocounter{fact}{7}
\setcounter{equation}{0}

%\end{document}

 \begin{center}

{\large\bf Addendum B. Analytic continuation of the Plancherel formula}

\end{center}

\renewcommand{\theequation}{B.\arabic{equation}}
\addtocounter{fact}{7}
\setcounter{equation}{0}

{\bf B.1. Plancherel formula.}
Plancherel formula for all classical series 1-10 of kernel-representations
 for large values of the parameter $\alpha$
there was obtained in \cite{Nerp} (this method is valid also for future
tubes).%
\footnote{For the groups $\U(p,q)$, $\Sp(2n,\R)$, $\SOS(2n)$, $\SO(n,2)$
and large values of parameter $\alpha$
the Plancherel formula was obtained by Berezin in 1978 (\cite{Ber2}),
Shortly after Berezin died
and proof was not published until \cite{UU}.
For the groups $\O(p,1)$, $\Sp(p,1)$ the Plancherel formula was obtained
in \cite{vD}}
 In all cases  the Plancherel measure has the form
\begin{equation}
V(\alpha) \prod_{k=1}^r \big| \Gamma(\frac12(\alpha -h+i s_k))\big|^2
 \frac {ds}{c(is)c(-is)} \mbox{\quad where\quad } \frac1i s_k\in\R
\end{equation}
Here $r$ is the rank of the group; $is_1$, $is_2$,...,$is_r$ are
the standard parameters of spherical functions,  spherical
functions of  unitary representations of non degenerate
principal  series corresponds to real $s_k$;  $\frac {ds} {c(is)c(-is)}$
is Gindikin--Karpelevich measure (see  \cite{Hel}), and
$V(\alpha)$ is a holomorphic factor.
The parameter $h$ is the least value
of $\alpha$ for which the kernel-representation is square integrable.

In this formula, the Plancherel measure is supported on
unitary principal series. We have seen that for small $\alpha$
spectrum of kernel-representation is complicated and hence  Plancherel
formula (B.1) is not correct for small $\alpha$.

{\bf B.2.
Analytic continuation of the Plancherel formula.}
Plancherel formula for $\alpha<h$ can be obtained by
the analytic continuation of formula (B.1).
We will give the final formula for the cases
$G=\O(p,q)$, $\U(p,q)$ (for other cases see \cite{Razlozhenie}).

Denote by  $g_t$ an element of $\O(p,q)$ having eigenvalues $e^{\pm t_k}$.
Let $\theta_0$ be the same as in Subsection 1.2.
Then
$$<T_\alpha(g_t) \theta_0, \theta_0>
=\prod_{k=1}^p { {\rm cosh}^{-\alpha}\, t_k}$$
Plancherel formula is equivalent to expansion of matrix element
 $<T_\alpha\theta_0,\theta_0>$ by spherical functions of $\O(p,q)$.
Let us denote spherical functions of $\O(p,q)$ by
 $\Phi_{is_1,\dots,is_r}$.  By  $\Phi_{is_1,\dots,is_r}(t_1,t_2,\dots t_p)$ we denote
a value of spherical function on $g_t$.

  In these notations the Plancherel formula for $\alpha>(p+q)/2-1$
 has the form
\begin{multline}
\prod_{k=1}^{p}  {\rm cosh}^{-\alpha} t_k= \\
=
   A\cdot\prod_{k=1}^p\frac{1}{\Gamma(\alpha-k+1)}
\int_{s\in\R^n}
\prod_{k=1}^p \Gamma ((\alpha-\frac12(p+q)+1+is_k)/2)
             \Gamma((\alpha-\frac12(p+q)+1-is_k)/2)
\times\\
\times
\prod_{k=1}^p \frac{\Gamma((q-p)/2+ is_k)\Gamma((q-p)/2-is_k)}
                      {\Gamma(is_k)\Gamma(-is_k)}
\times \\
\prod_{1\le k< l\le p} (s_k^2-s_l^2)\tanh \frac{\pi}{2}(s_k-s_l)
                                       \tanh \frac{\pi}{2}(s_k+s_l)
\times \\
\times
\Phi_{is_1,\dots,is_r}(t_1,t_2,\dots, t_r)\, ds_1 \dots ds_p %\notag
\end{multline}
where $A$ is a constant.

\medskip

Now we will wright the analytic continuation of formula
(B.2). Let $r=0,1,\dots, p$.
Let $u_1,u_2, \dots,u_r$ be nonnegative integers.
Let
$$w_j=u_1+\dots+u_j+j/2$$

Then
\begin{multline}
\prod_{m=1}^p {\rm cosh}^{-\alpha} (t_m)=\\ =
\sum_{r=0}^p\qquad \sum_{u_1,\dots,u_r: w_r<h-\alpha}
C_u\cdot V_u(\alpha) \times \\ \times
\int Q_u(\alpha|s_{r+1},\dots,s_p)
\Phi_{\alpha-(p+q)/2+2w_1,\dots,\alpha-(p+q)/2+2w_r,is_{r+1},\dots,is_p}
(t_1,\dots,t_p)\times \\
ds_{r+1}\dots ds_{p}
\end{multline}
where
$$
C=A\frac{2^{p-r}p! (2\pi)^r}{(p-r)!}
\frac{(-1)^{u_1+\dots+u_s}}{\prod u_j!}\cdot
\prod_{k<m\le s} \frac{\Gamma(\frac12+w_m-w_k)} {\Gamma(w_m-w_k)
  (\frac12+w_{k-1}-w_m)_{u_k}}
$$
and $A$ is a constant,
\begin{multline*}
V_u(\alpha)= \frac{1}{\prod_{m=1}^p\Gamma(\alpha-k+1)} \times \\
\prod_{k=1}^r\frac{\Gamma(\alpha-p+1++2w_k)\Gamma(-\alpha+q-1-2w_k)}
                  {\Gamma(-\alpha+\frac12(p+q)-2w_k)
                  (\alpha-\frac12(p+q)+w_k+w_{k-1}+\frac12)_{u_k}}\times \\
\prod_{k<m\le r} \frac{\Gamma(\frac12-\alpha+\frac12(p+q)-w_k-w_m)}
                  {\Gamma(-\alpha+\frac12(p+q)-w_k-w_m)
                   (\alpha-\frac12(p+q)+w_m+w_{k-1}+\frac12)_{u_k}}
\end{multline*}
and
\begin{multline*}
Q_u(\alpha|s_{r+1},\dots,s_p)=\\=
\prod_{n=r+1}^p |\Gamma(\frac12(\alpha-(p+q)/2+1+is_n)+w_{n-1})|^2\times \\
\prod_{k\le r, l>r}\left|
      \frac{\Gamma(\frac12(1-\alpha+\frac12(p+q)-2w_k+is_l)}
       {\Gamma(\frac12(-\alpha+\frac12(p+q)-2w_k+is_l))
        (\frac12(\alpha-\frac12(p+q)+w_k +is_l)_{u_k}}
     \right|^2
\times\\
\prod_{l>n>r} (s_l^2-s_n^2) \tanh(\frac{\pi}2(s_l-s_n))\tanh(\frac {\pi}2(s_l+s_n))
\end{multline*}

In the case $G=\U(p,q)$ we obtain
\begin{multline}
\prod_{m=1}^p {\rm cosh}^{-\alpha} (t_m)=\\ =
\sum_{r=0}^p\qquad \sum_{w_1,\dots,w_r}
C_w\cdot V_w(\alpha) \times \\ \times
\int Q_w(\alpha|s_{r+1},\dots,s_p)
\Phi_{\alpha-(p+q)+1+2w_1,\dots,\alpha-(p+q)+1+2w_r,is_{r+1},\dots,is_p}
(t_1,\dots,t_p)\times \\
ds_{r+1}\dots ds_{p}
\end{multline}
where
$$
C_w=A\cdot\frac{2^{p-r}p! (2\pi)^r}{(p-r)!}\cdot
\prod_{k=1}^r\frac{(-1)^{w_k}}{w_k!}
\prod_{k<l\le r}(w_k-w_l)^2
$$
$A$ is a constant,
\begin{multline*}
V_w(\alpha)=\\=
\prod_{m=1}^p\frac{1}{\Gamma(\alpha/2-m+1)^2}
\prod_{k=1}^r\frac{\Gamma(\frac12\alpha-p+1+w_k)^2
 \Gamma(-\frac12\alpha+ q-w_k)^2}
{\Gamma(-\alpha+(p+q)-1+2w_k)(\alpha-(p+q)+1+w_k)_{w_k}}\times\\
\prod_{k<l\le r}(\alpha-(p+q)+1+w_k +w_l)^2
\end{multline*}
and
\begin{multline*}
Q_w(\alpha|s_{r+1},\dots,s_p)=\\=
\prod_{k\le r, n>r}
 ( \frac14(\alpha-(p+q)+1+2w_k)^2+s_n^2)^2 \times\\
\prod_{n=r+1}^p
 \frac{|\Gamma(\frac12((q-p+1+is_n))|^4 |\Gamma(\frac12(\alpha- p-q+1+is_n)|^2}
                    {|\Gamma(is_n)|^2}   \times\\
\prod_{n>m>r} (s_n^2-s_m^2)^2
\end{multline*}

{\bf B.3. Negative integer $\alpha$.} Analytic formula (B.3),(B.4)
 are valid for all values of the parameter
$\alpha\in\R$. Nevertheless it has visible sense only for
$$ \alpha=0,1,2,\dots,p-1 \mbox{\quad or \quad} \alpha>p-1$$
when representation $T_\alpha$ is unitary, and for
$$\alpha=0,-1,-2,\dots$$
 when the representation $T_\alpha$ is finite-dimensional.
In the last case  the right parts of identities contain only finite number
of summands with $r=p$ ( $\Gamma$-factors
in denominator annihilate all summands with $r<p$).
In the case $G=\U(p,p)$ we obtain Pickrell decomposition
, \cite{Pic}.
In this place our work also has an intersection with a  recent preprint of
G.Zhang \cite{Zhang}.

 For over values of parameter $\alpha$ Plancherel formula (B.3), (B.4)
 are
 elements of mysterious  nonunitary harmonic analysis
(see, for instance \cite{Molf},\cite{MD}).

\renewcommand{\theequation}{C.\arabic{equation}}
\addtocounter{fact}{7}
\setcounter{equation}{0}

\end{document}